\documentclass[12pt, twoside, leqno]{article}

\usepackage{amsmath,amsthm}
\usepackage{amssymb}

\usepackage{enumitem}

\usepackage{graphicx}

\usepackage[T1]{fontenc}

\usepackage{amsmath,amstext,amsthm,amscd,typearea}
\usepackage{amssymb}
\usepackage{a4wide}
\usepackage[mathscr]{eucal}
\usepackage{mathrsfs}
\usepackage{typearea}
\usepackage{charter}
\usepackage{pdfsync}
\usepackage{url}

\usepackage{xcolor}

\usepackage[a4paper,width=16.2cm,top=3cm,bottom=3cm]{geometry}

\numberwithin{equation}{section}

\newtheorem{theorem}{Theorem}[section]

\newtheorem{proposition}[theorem]{Proposition}
\newtheorem{corollary}[theorem]{Corollary}
\newtheorem{lemma}[theorem]{Lemma}
\newtheorem{remark}[theorem]{Remark}

\newcommand{\cali}[1]{\mathscr{#1}}

\newcommand{\ddc}{dd^c}
\newcommand{\dc}{d^c}

\newcommand{\C}{\mathbb{C}}

\newcommand{\N}{\mathbb{N}}

\newcommand{\R}{\mathbb{R}}

\begin{document}


\baselineskip=17pt


\title{Convexity of the class of currents with finite relative energy}

\author{Duc-Viet Vu,\\ University of Cologne,\\ Mathematics Division,\\ Department of Mathematics and Computer Science,\\
 Weyertal 86-90, 50931 K\"oln,  Germany\\
E-mail: vuviet@math.uni-koeln.de}

\date{}

\maketitle


\renewcommand{\thefootnote}{}

\footnote{2020 \emph{Mathematics Subject Classification}: Primary 32U15.}

\footnote{\emph{Key words and phrases}: relative non-pluripolar product, relative energy, full mass intersection.}

\renewcommand{\thefootnote}{\arabic{footnote}}
\setcounter{footnote}{0}


\begin{abstract}
We prove  the convexity of the class of currents with finite relative energy.   A key ingredient is an integration by parts formula for  relative non-pluripolar products which is of independent interest.
\end{abstract}


\section{Introduction}


 Let $X$ be a compact K\"ahler manifold of dimension $n$.  Let $T$ be a closed positive current of bi-degree $(p,p)$ on $X$. Let $T_1, \ldots, T_m$ be closed positive $(1,1)$-currents on $X$.   \emph{The $T$-relative non-pluripolar product} $\langle T_1 \wedge \cdots \wedge T_m \dot{\wedge} T \rangle$ of $T_1,\ldots, T_m$  was introduced in \cite{Viet-generalized-nonpluri}. The last product is a closed positive current of bi-degree $(p+m,p+m)$.    When $T$ is a constant function equal to $1$ (\emph{i.e,} $T$ is the current of integration along $X$), the current $\langle T_1 \wedge \cdots \wedge T_m \dot{\wedge} T \rangle$  coincides with the usual non-pluripolar product of $T_1, \ldots, T_m$  given in \cite{BT_fine_87,BEGZ,GZ-weighted}. 
 
For every closed positive currents $S$ on $X$, we denote by $\{S\}$ its cohomology class.  For two cohomology $(q,q)$-classes $\alpha,\beta$ on $X$, we write $\alpha \le \beta$ if $\beta- \alpha$ can be represented by a closed positive $(q,q)$-current.
   
Recall that by \cite[Theorem 1.1]{Viet-generalized-nonpluri} (also \cite{BEGZ,WittNystrom-mono,Lu-Darvas-DiNezza-mono}), if $T'_j$ is a closed positive  $(1,1)$-current on $X$ which is cohomologous to $T_j$ and  less singular than $T_j$ for $1 \le j \le m$, then we have 
\begin{align}\label{ine-mono-relativenopluri}
\{\langle T_1 \wedge \cdots \wedge T_m \dot{\wedge} T \rangle \} \le \{\langle T'_1 \wedge \cdots \wedge T'_m \dot{\wedge} T \rangle \}.
\end{align}
The last inequality allows us to define the notion of $T$-relative full mass intersection, see \cite{Viet-generalized-nonpluri,BEGZ}.  We say that $T_1,\ldots, T_m$ are of \emph{$T$-relative full mass intersection} if 
$$\{\langle \bigwedge_{j=1}^m T_j \dot{\wedge} T\rangle\}=\{\langle \bigwedge_{j=1}^m T_{j,\min} \dot{\wedge} T\rangle\},$$
 where $T_{j, \min}$ is a current with minimal singularities in the class $\{T_j\}$ for $1 \le j \le m$. 

Let $\alpha$ be a pseudoeffective $(1,1)$-class.  Denote by  $\mathcal{E}_{m}(\alpha, T)$ the set of currents $P \in \alpha$ such that $P, \ldots, P$ ($m$ times $P$) are of $T$-relative full mass intersection. 

Recall that $\mathcal{W}^-$ is  the set of convex increasing functions $\chi$ from $\R$ to $\R$ such that $\chi(-\infty)= -\infty$. Let $\chi \in \mathcal{W}^-$.   We can  define $\mathcal{E}_{\chi,m}(\alpha,T)$ to be the subclass of $\mathcal{E}_{m}(\alpha, T)$ consisting of $P$ such that $P$ has finite  ($T$-relative) $\chi$-energy.  When $T \equiv 1$ and $m=n$, the class $\mathcal{E}_{\chi,m}(\alpha,T)$ generalizes the usual class of currents with finite energy in \cite{BEGZ,GZ-weighted}, see also \cite{Cegrell} for the local setting. We refer to Section \ref{sec-weighted} for details. For the moment,  we note here that 
$$\mathcal{E}_{m}(\alpha,T)= \bigcup_{\chi \in \mathcal{W}^-}\mathcal{E}_{\chi,m}(\alpha,T).$$     
 Here is our main result. 

\begin{theorem}\label{th-convexrelativeclassintro} The sets $\mathcal{E}_{\chi,m}(\alpha, T)$ and  $\mathcal{E}_{m}(\alpha,T)$ are convex.
\end{theorem}

The last result was proved in \cite[Theorem 1.3]{Viet-generalized-nonpluri} in the case where $\alpha$ is K\"ahler.   When $m=n$ and  $T\equiv 1$, the convexity  $\mathcal{E}_{\chi,m}(\alpha, T)$ was conjectured in \cite{BEGZ}. It was later answered affirmatively in \cite[Corollary 2.12]{Lu-Darvas-DiNezza-singularitytype} in this setting. The proof in \cite{Lu-Darvas-DiNezza-singularitytype} doesn't extend directly to our setting because it uses, in a crucial way,  Monge-Amp\`ere equations in big classes. 


We will see that Theorem \ref{th-convexrelativeclassintro} is a direct consequence of  a more general result (Theorem \ref{th-convexrelativeclass}) which is in turn deduced from  a monotonicity property of joint energy of currents, see Theorem \ref{th-jointenergy} below. To prove these results, we  use ideas from  the proof of \cite[Theorem 1.3]{Viet-generalized-nonpluri} and prove an integration by parts formula for relative non-pluripolar products (Theorem \ref{th-integra}) which is of independent interest.   \emph{We emphasize that the last formula was applied to the study of complex Monge-Amp\`ere equations.} It plays a key role in the proof of main results in \cite{Vu_Do-MA}, see Theorem 1.3 there. 

Moreover it was also explained in \cite{Vu_Do-MA} that by using the integration by parts formula obtained in this work and the variational method (\cite{BBGZ-variational,Lu-Darvas-DiNezza-mono}), one can solve the Monge-Amp\`ere equation in the prescribed singularity setting without the small unbounded locus assumption. Hence this gives another proof of a main result in \cite{Lu-Darvas-DiNezza-logconcave}.  We refer to \cite[Theorem 3.8]{Vu_Do-MA} for details.  

In the next section, we will present the above-mentioned integration by parts formula for relative non-pluripolar products. This formula strengthens (and generalizes) recent ones obtained in \cite{Lu-comparison-capacity,Xia} (see Corollary \ref{cor-integraintro} and the paragraph following it). Our main result will be proved in Section \ref{sec-weighted}. \\


\section{Integration by parts}

We first recall some basic facts about relative non-pluripolar products. This notion was introduced in \cite{Viet-generalized-nonpluri} as a generalization of the usual non-pluripolar products  given in \cite{BT_fine_87,BEGZ,GZ-weighted}.  

Let $X$ be a compact K\"ahler manifold.  Let $T_1, \ldots, T_m$ be closed positive $(1,1)$-currents on $X$. Let $T$ be  a closed positive current of bi-degree $(p,p)$ on $X$. By \cite{Viet-generalized-nonpluri}, the $T$-relative non-pluripolar product $\langle \bigwedge_{j=1}^m T_j \dot{\wedge} T\rangle$ is defined  in a way similar to that of  the usual non-pluripolar product. The product $ \langle  \bigwedge_{j=1}^m T_j \dot{\wedge} T\rangle $ is a well-defined closed positive current of bi-degree $(m+p,m+p)$; and $ \langle  \bigwedge_{j=1}^m T_j \dot{\wedge} T\rangle $  is  symmetric with respect to $T_1, \ldots, T_m$ and is homogeneous.

For every closed positive $(1,1)$-current $P$, we denote by $I_P$ the set of $x \in X$ so that local potentials of $P$ are equal to $-\infty$ at $x$. Note that $I_P$ is  a locally complete pluripolar set. The following is deduced from \cite[Proposition 3.5]{Viet-generalized-nonpluri}.

\begin{proposition}\label{pro-sublinearnonpluripolar}

$(i)$ Given a locally complete pluripolar set $A$ such that $T$ has no mass on $A$, then $\langle \bigwedge_{j=1}^m T_j \dot{\wedge} T\rangle$ also has no mass on $A$.

$(ii)$  Let $T'_1$ be  a closed positive $(1,1)$-current on $X$ and $T_j,T$ as above. Assume that  $T$ has no mass on $I_{T_1} \cup I_{T'_1}$.  Then we have
\begin{align}\label{ine-convexnonpluripoarTjTphayj}
\big\langle  (T_1+T'_1) \wedge \bigwedge_{j=2}^m T_j \dot{\wedge} T\big\rangle = \langle T_1 \wedge \bigwedge_{j=2}^m T_j \dot{\wedge}T\rangle+ \langle T'_1 \wedge  \bigwedge_{j=2}^m  T_j \dot{\wedge} T\rangle.
\end{align}

$(iii)$ Let $1 \le l \le m$ be an integer. Then for $R:= \langle \bigwedge_{j=l+1}^m T_j \dot{\wedge} T \rangle$, there holds $\langle \bigwedge_{j=1}^m T_j \dot{\wedge} T \rangle = \langle \bigwedge_{j=1}^l T_j \dot{\wedge} R \rangle$.

$(iv)$ The equality
$$\langle \bigwedge_{j=1}^m T_j \dot{\wedge} T \rangle = \langle \bigwedge_{j=1}^m T_j \dot{\wedge} T' \rangle$$
holds, where $T':= \bold{1}_{X \backslash \bigcup_{j=1}^m I_{T_j}} T$.
\end{proposition}

We also need  the following result. 

\begin{theorem} \label{the-th-increasingsequenceMa} (\cite[Theorem 2.6, Remark 2.7]{Viet-generalized-nonpluri}) Let $u_j$ be a locally bounded plurisubharmonic (psh) function on  an open subset $U$ of $\C^n$ for $1 \le j \le m$. Let $(u_{jk})_{k \in \N}$ be a sequence of locally bounded psh functions increasing  to $u_j$  almost everywhere as $k \to \infty$. Let $T$ be a closed positive current on $U$.  Then, the convergence 
\begin{align}\label{limit-increseinpointwise}
u_{1k} \ddc u_{2k} \wedge \cdots \wedge \ddc u_{mk} \wedge T \to u_1 \ddc u_{2} \wedge \cdots \wedge \ddc u_{m} \wedge T
\end{align}
as $k \to \infty$ holds provided that $T$ has no mass on $A_j:= \{x \in U:  u_j(x) \not = \lim_{k\to \infty} u_{jk}(x) \}$ for every $1 \le j \le m$ and  the set $A_j$ is locally complete pluripolar for every $j$. 
\end{theorem}

Recall that a \emph{dsh} function on $X$ is the difference of two quasiplurisubharmonic (quasi-psh for short) functions on $X$ (see \cite{DS_tm}). These functions are well-defined outside pluripolar sets. Let $v$ be a dsh function on $X$. Write $v= \varphi_1- \varphi_2$, where $\varphi_1, \varphi_2$ are quasi-psh function. Hence $v(x)$ is well-defined  for  $x \in X \backslash A$, where $A:= \{\varphi_1= -\infty\}\cup \{\varphi_2= -\infty\}$.   The function $v$ is said to be \emph{bounded} in $X$ if there exists a constant $C$ such that $|v(x)| \le C$ for every $x \in X \backslash A$.  

We say that $v$ is \emph{$T$-admissible} (or admissible with respect to $T$) if  there exist  quasi-psh functions $\varphi_1, \varphi_2$ on $X$ such that $v= \varphi_1- \varphi_2$  and $T$ has no mass on $\{\varphi_j=-\infty\}$ 
for $j=1,2$. In particular, if $T$ has no mass on pluripolar sets, then every dsh function is $T$-admissible. Assume now that $v$ is $T$-admissible.  The following is a direct consequence of  Proposition \ref{pro-sublinearnonpluripolar} $(i)$.

\begin{lemma}\label{le-phiadmissible} If $v$ is $T$-admissible, then $v$ is  also admissible with respect to $\langle \bigwedge_{j=1}^m T_j \dot{\wedge} T \rangle$. 
\end{lemma}

Recall that if  $v= \varphi_1- \varphi_2$ for some bounded quasi-psh functions $\varphi_1, \varphi_2$ on $X$ (note $X$ is compact), the current $dv \wedge \dc v \wedge T$ is, by definition, equal to 
\begin{align}\label{eq-dinhnghiadvdcv}
\frac{1}{2}\ddc (\varphi_1- \varphi_2)^2  \wedge T - (\varphi_1- \varphi_2) \ddc (\varphi_1- \varphi_2) \wedge T.
\end{align}
We notice that in the above formula the function $(\varphi_1- \varphi_2)^2$ is locally the difference of two bounded psh functions. More precisely, we can assume $\varphi_j$ is $\omega$-psh function for $j=1,2$, where $\omega$ is a K\"ahler form on $X$, and   consider an open local chart $U$ such that $\omega = \ddc \phi$ for some smooth psh function $\phi$ on $U$. By adding to $\phi$ a big constant, we can even assume that $\varphi_j+ \phi \ge 0$ on $U$ for $j=1,2$. Thus $(\varphi_j+\phi)^2$ and $(\varphi_1+ \varphi_2+2 \phi)^2$ are psh on $U$ for $j=1,2$, and 
\begin{align}\label{eq-bieudienhieuquasi}
(\varphi_1- \varphi_2)^2= (2 (\varphi_1+ \phi)^2 + 2 (\varphi_2+\phi)^2)- (\varphi_1+ \varphi_2+2 \phi)^2
\end{align}
which is the difference of two bounded psh functions on $U$. 

Consider another dsh function $w$ which is equal to the difference of two locally bounded psh functions and $T$ is of bi-degree $(n-1, n-1)$, we have 
\begin{align}\label{eq-tichdvdcw}
2 d v \wedge \dc w \wedge T= d(v+ w) \wedge \dc (v+ w) \wedge T - d v \wedge \dc v\wedge T- d w \wedge \dc w \wedge T.
\end{align}
However, in general, even when $v$ is bounded, $v$ might not be the difference of two bounded quasi-psh functions.  Hence, the current $``dv \wedge \dc v \wedge T "$ is not well-defined in the above sense. We will introduce below the current $\langle dv\wedge \dc v \dot{\wedge} T \rangle$ in the spirit of non-pluripolar products.   Before going into details, we need the following auxiliary estimate.
 
\begin{lemma}\label{le-dvdcv} Let $\omega$ be a K\"ahler form on $X$. Let $\varphi_1,\varphi_2$ be bounded $\omega$-psh functions on $X$ and $v:= \varphi_1- \varphi_2$.  Let $T$ be a closed positive current of bi-dimension $(1,1)$ on $X$. Then, there exists a constant $C$ independent of $\varphi_1,\varphi_2$ such that
\begin{align}\label{ine-danhgiadvdcv}
\int_X d v \wedge \dc v \wedge T \le C \|v\|_{L^\infty}.
\end{align}
\end{lemma}

\proof  We have 
\begin{align*}
I &:= \int_X  dv \wedge \dc v \wedge T= -  \int_X v \ddc v \wedge T\\
& = -\int_X v (\ddc \varphi_1 -\ddc \varphi_2) \wedge T\\
& = -\int_X v (\ddc \varphi_1+ \omega)\wedge T+ \int_X v (\ddc \varphi_2+ \omega)\wedge T\\
& \le \|v\|_{L^\infty}\sum_{j=1}^2 \int_X (\ddc \varphi_j+ \omega)\wedge T\\
&= 2\|v\|_{L^\infty} \int_X T \wedge \omega
\end{align*}
by Stokes' theorem.  The desired estimate follows.  The proof is finished.
\endproof

Assume now that $v$ is $T$-admissible.    Let $\varphi_{1}, \varphi_{2}$ be quasi-psh functions such that $v= \varphi_{1}- \varphi_{2}$ and $T$ has no mass on $\{\varphi_{j}=-\infty\}$ for $j=1,2$. Let $\varphi_{j,k}:= \max\{\varphi_{j}, -k \}$ for every $j=1,2$ and $k \in \N$. Put $v_k:= \varphi_{1,k}- \varphi_{2,k}$.   Since $v_k$ is the difference of two bounded quasi-psh functions, using (\ref{eq-dinhnghiadvdcv}),   we obtain
$$Q_k:= d v_k \wedge \dc v_k \wedge T=\ddc v_k^2  \wedge T - v_k\ddc v_k \wedge T.$$
Let $\omega$ be a K\"ahler form so that $\varphi_j$ is $\omega$-psh for $j=1,2$. Let $U$ be a local chart on $X$ such that $\omega= \ddc \phi$ on $U$ for some psh function $\phi$ such that $\varphi_{j,k}+ \phi \ge 0$ on $U$ for $j=1,2$ (we fix $k$).  By (\ref{eq-bieudienhieuquasi}) applied to $\varphi_{j,k}+\phi$, we have
\begin{align*}
Q_k &= 2 \ddc (\varphi_{1,k}+\phi)^2 \wedge T+2 \ddc (\varphi_{2,k}+\phi)^2 \wedge T - \ddc(\varphi_{1,k}+ \varphi_{2,k}+2\phi)^2 \wedge T\\ & \quad - 
(v_k \ddc \varphi_{1,k} \wedge T-v_k \ddc \varphi_{2,k} \wedge T)
\end{align*}
on $U$. By the plurifine locality with respect to $T$ (\cite[Theorem 2.9]{Viet-generalized-nonpluri}) applied to each term in the right-hand side of the last equality, we have 
\begin{align}\label{eq-localplurifineddc}
\bold{1}_{\bigcap_{j=1}^2 \{\varphi_{j}> -k\}} Q_k =\bold{1}_{\bigcap_{j=1}^2\{\varphi_{j}> -k\}} Q_{k'}
\end{align}
for every $k'\ge k$.

We say that  $\langle d v \wedge \dc v \dot{\wedge} T \rangle$ is \emph{well-defined} if   the mass of $\bold{1}_{\bigcap_{j=1}^2 \{\varphi_{j}> -k\}} Q_k$ on $X$ is bounded uniformly in $k$. In this case,  using  (\ref{eq-localplurifineddc}) implies that there exists a positive current $Q$ on $X$ such that for every bounded Borel form $\Phi$ on $X$, we have
$$\langle Q,\Phi \rangle  = \lim_{k\to \infty} \langle \bold{1}_{\bigcap_{j=1}^2 \{\varphi_{j}> -k\}} Q_k, \Phi\rangle.$$
We define $\langle d v \wedge \dc v \dot{\wedge} T \rangle$ to be the current $Q$. This agrees with the classical definition if $v$ is the difference of two bounded quasi-psh functions. This definition is independent of the choice of $\varphi_1, \varphi_2$ by Lemma \ref{le-duynhat} below. If $w$ is another $T$-admissible dsh function and $T$ is of bi-dimension $(1,1)$ such that the currents $\langle d v \wedge \dc v \dot{\wedge} T \rangle$, $\langle d w \wedge \dc w \dot{\wedge} T \rangle$, and  $\langle d (v+w) \wedge \dc (v+w) \dot{\wedge} T \rangle$ are all well-defined, we define $\langle d v \wedge \dc w \dot{\wedge} T \rangle$ using (\ref{eq-tichdvdcw}) formally. 

\begin{lemma} \label{le-duynhat} Let $\varphi'_1, \varphi'_2$ be quasi-psh functions on $X$ such that $v= \varphi'_1- \varphi'_2$ and $T$ has no mass on $\{\varphi_j'= - \infty\}$ for $j=1,2$. Let $\varphi'_{j,k}, Q'_k$ be  the function and  current associated to $\varphi'_j$ defined similarly as $\varphi_{j,k}$ and $Q_k$ respectively. Then if $\bold{1}_{\bigcap_{j=1}^2 \{\varphi_{j}> -k\}} Q_k$ is of mass bounded uniformly on $k$ then so is $\bold{1}_{\bigcap_{j=1}^2 \{\varphi'_{j}> -k\}} Q'_k$, and 
\begin{align}\label{eq-Qduynhat}
Q= \lim_{k\to \infty} \langle \bold{1}_{\bigcap_{j=1}^2 \{\varphi'_{j}> -k\}} Q'_k, \Phi\rangle
\end{align}
for every bounded Borel form $\Phi$ on $X$.
\end{lemma} 

\proof Since $v=\varphi'_1- \varphi_2'$, we get 
$$\varphi_1+ \varphi'_2= \varphi'_1+ \varphi_2.$$
Put $v'_k:= \varphi'_{1,k}- \varphi'_{2,k}$,   $A'_k:= \bigcap_{j=1}^2 \{\varphi'_{j}> -k\}$, and $A_k:= \bigcap_{j=1}^2 \{\varphi_{j}> -k\}$. We have $Q'_k= d v'_k \wedge \dc v'_k \wedge T$, and  $v_k= v'_k$ on $A_k \cap A'_k$ which is open in the plurifine topology. We claim that 
\begin{align}\label{eq-QkbangnQkphay}
\bold{1}_{A_k \cap A'_k}Q'_k =  \bold{1}_{A_k \cap A'_k} Q_k.
\end{align}
This is a sort of plurifine locality statement and can be essentially deduced from the plurifine locality for bounded psh functions (here we have $v_k= v'_k$ on $A_k \cap A'_k$ but $v_k, v'_k$ are only dsh). We give details for readers' convenience. Before doing so, we will show that the desired assertion is a direct consequence of (\ref{eq-QkbangnQkphay}). First observe that $\bold{1}_{A'_k}Q'_k$ has no mass on the pluripolar set $\{\varphi_j= -\infty\}$ for $j=1,2$ by Proposition \ref{pro-sublinearnonpluripolar} $(i)$ and the fact that $T$ has no mass on $\{\varphi_j= -\infty\}$. It follows that    
\begin{align*}
\bold{1}_{A'_k}Q'_k &= \lim_{s \to \infty} \bold{1}_{A'_k \cap A_{s}} Q'_k \\
&=\lim_{s \to \infty} \bold{1}_{A'_k \cap A_{s}} Q'_s\\
& \le \lim_{s \to \infty} \bold{1}_{A'_s \cap A_{s}} Q'_s\\
& =\lim_{s \to \infty} \bold{1}_{A'_s \cap A_{s}} Q_s  \le  \lim_{s\to\infty} \bold{1}_{A_{s}} Q_s= Q,
\end{align*}
where we used (\ref{eq-localplurifineddc}) applied to $Q'_k$ in the second equality and  (\ref{eq-QkbangnQkphay}) in the third equality. By exchanging the role of $Q'_k$ and $Q_k$, we also obtain that 
$$\bold{1}_{A_k} Q_k \le R,$$ 
for every limit current $R$ of $(\bold{1}_{A'_k}Q'_k)_k$ as $k \to \infty$. Hence (\ref{eq-Qduynhat}) follows.

We go back to the proof of (\ref{eq-QkbangnQkphay}). Write
$$d v_k \wedge \dc v_k \wedge T= d v_k \wedge \dc (v_k - v'_k) \wedge T + d(v_k- v'_k) \wedge \dc v'_k \wedge T+ d v'_k \wedge \dc v'_k \wedge T.$$
Denote by $R_1,R_2$ the first and second currents in the right-hand side of the last equality. In order to obtain (\ref{eq-QkbangnQkphay}), it suffices to check that $R_j= 0$ on $A_k \cap A'_k$ for $j=1,2$. Observe
\begin{align*}
R_1 &= d v_k \wedge \dc (\varphi_{1,k}+ \varphi_{2,k}' - \varphi_{1,k}' - \varphi_{2,k}) \wedge T \\
& = \big[d \varphi_{1,k} \wedge \dc (\varphi_{1,k}+ \varphi_{2,k}')\wedge T -  d \varphi_{1,k} \wedge \dc (\varphi'_{1,k}+ \varphi_{2,k})\wedge T \big] -\\
& \quad \big[ d \varphi_{2,k} \wedge \dc (\varphi_{1,k}+ \varphi_{2,k}')\wedge T -  d \varphi_{2,k} \wedge \dc (\varphi'_{1,k}+ \varphi_{2,k})\wedge T \big].
\end{align*}   
Each term in the right-hand side of the above equality is equal to $0$ on $A_k \cap A'_k$ thanks to the plurifine locality and the fact that
$$\varphi_{1,k}+ \varphi'_{2,k}= \varphi_1+ \varphi'_2= \varphi_1'+ \varphi_2=\varphi_{1,k}'+ \varphi_{2,k}$$
on $A_k \cap A'_k$. Hence $R_1= 0$ on $A_k \cap A'_k$. Similarly we get $R_2= 0$ on $A_k \cap A'_k$. This finishes the proof.
\endproof

\begin{lemma} \label{le-existenceddcT} Assume that $v$ is bounded. Then, the current $\langle d v \wedge \dc v \dot{\wedge} T \rangle$ is well-defined.
\end{lemma}

\proof  Let the notation be as in the proof of Lemma \ref{le-dvdcv}. Let $\omega$ be a K\"ahler form on $X$ such that $\varphi_1, \varphi_2$ are $\omega$-psh. Note that $\varphi_{jk}$ is also $\omega$-psh for every $j,k$.   Observe that since $v$ is bounded,  there exists a constant $C$ such that $\varphi_2-C \le \varphi_1 \le \varphi_2+ C$. Thus,  there exists a constant $C$ so that 
\begin{align*} 
\|v_k\|_{L^\infty} \le C
\end{align*}
for every $k$.  Using this and Lemma \ref{le-dvdcv},  one gets
$$\|Q_k\| \lesssim \|v_k\|_{L^\infty}  \le C$$
for some constant $C$ independent of $k$. Hence, the desired assertion follows. This finishes the proof.
\endproof

Let $T_1, \ldots, T_m$ be  closed positive $(1,1)$-currents on $X$ and $R:= \langle T_1 \wedge \cdots \wedge T_m \dot{\wedge} T\rangle$. We define 
$$\langle  dv \wedge \dc v \wedge T_1 \wedge \cdots \wedge T_m \dot{\wedge} T\rangle:= \langle d v \wedge  \dc v \dot{\wedge} R\rangle.$$
When $T \equiv 1$, we write the left-hand side of the last equality simply as $\langle  dv \wedge \dc v \wedge T_1 \wedge \cdots \wedge T_m\rangle$.

The current $\langle  dv \wedge \dc w \wedge T_1 \wedge \cdots \wedge T_m \dot{\wedge} T\rangle$
is defined similarly if $p+m=n-1$, where $T$ is of bi-degree $(p,p)$.  We put
$$\langle \ddc v \dot{\wedge} T \rangle:= \langle \ddc \varphi_1 \dot{\wedge} T\rangle- \langle \ddc \varphi_2 \dot{\wedge} T\rangle.$$
 Define 
$$\langle  \ddc v\wedge T_1 \wedge \cdots \wedge T_m \dot{\wedge} T\rangle:= \langle  \ddc v\dot{\wedge}R \rangle.$$

By Proposition \ref{pro-sublinearnonpluripolar} $(iii)$, this definition agrees with the $T$-relative non-pluripolar product of $\ddc v, T_1, \ldots, T_m$ if $v$ is quasi-psh. When $T \equiv 1$, we write $\langle  \ddc v\wedge T_1 \wedge \cdots \wedge T_m \rangle$ for $\langle  \ddc v\wedge T_1 \wedge \cdots \wedge T_m \dot{\wedge} T\rangle$. In this case the product $\langle  \ddc v\wedge T_1 \wedge \cdots \wedge T_m \rangle$ is the one defined in the paragraph right after Theorem 1.2 in \cite{Lu-comparison-capacity}.

By admissibility and Proposition \ref{pro-sublinearnonpluripolar} $(ii)$, we can check that if $v, w$ are dsh functions which are admissible with respect to $T$, then 
$$\langle  \ddc (v+w)  \dot{\wedge} T\rangle=\langle  \ddc v\dot{\wedge} T\rangle+\langle  \ddc w   \dot{\wedge} T\rangle.$$ 

Here is an integration by parts formula for relative non-pluripolar products. 

\begin{theorem} \label{th-integra} Let  $T$ be a closed positive current of bi-degree $(n-1,n-1)$ on $X$. Let $v,w$ be bounded $T$-admissible dsh functions on $X$. Then,  we have  
\begin{align}\label{eq-intebyparts}
\int_X w \langle  \ddc v \dot{\wedge} T\rangle=\int_X v \langle \ddc w\dot{\wedge} T\rangle=- \int_X \langle dw \wedge \dc v \dot{\wedge} T \rangle.
\end{align}
\end{theorem}

The last result was proved in \cite[Theorem 1.14]{BEGZ} if $v,w$ can be written as the differences of psh functions which are locally bounded outside a \emph{closed} locally complete pluripolar set; see also \cite{Bedford_Taylor_82,Sibony_duke}.  

\begin{proof} We use ideas from the proof of \cite[Proposition 4.2]{Viet-generalized-nonpluri}. 
 Let $\varphi_1,\varphi_2, \varphi_3,\varphi_4$ be  negative quasi-psh functions on $X$ such that  $v= \varphi_1- \varphi_2$ and $w= \varphi_3- \varphi_4$ and $T$ has no mass on $\bigcup_{j=1}^4 \{\varphi_j = - \infty\}$. Let $\omega$ be a K\"ahler form on $X$ such that $\varphi_j$ is $\omega$-psh for every $1 \le j \le 4$.   Put
 $$\psi:= \varphi_1+\varphi_2+\varphi_3+\varphi_4, \quad \psi_k:= k^{-1}\max\{\psi, -k\}+1.$$
 and $\varphi_{jk}:= \max\{\varphi_j, -k\}$ for $1 \le j \le 4$. Observe that $0 \le \psi_k \le 1$. Let  $x \in X$ such that $\psi_k(x) >0$. We have 
 $$\varphi_1(x)+\varphi_2(x)+\varphi_3(x)+\varphi_4(x)= \psi(x) >-k.$$
This combined with the property that $\varphi_j \le 0$ for every $1 \le j \le 4$ yields that $\varphi_j(x)>-k$ for every $1 \le j \le 4$. We infer that 
\begin{align}\label{inclu-psikkhacko}
\{\psi_k  \not = 0\} \subset \bigcap_{j=1}^4\{\varphi_j >-k\}.
\end{align}

Put $v_k:= \varphi_{1k}- \varphi_{2k}$ and $w_k:=  \varphi_{3k}- \varphi_{4k}$. Since $v$ and $w$ are bounded, the functions $v_k,w_k$ are bounded uniformly in $k$.

Let $A:=\bigcup_{j=1}^4\{\varphi_j=-\infty\}$.  By admissibility and Proposition \ref{pro-sublinearnonpluripolar} $(i)$, we see that  
\begin{align} \label{eq-masstrenpluripolar}
\bold{1}_{A}\langle (\ddc \varphi_j+\omega) \dot{\wedge} T\rangle=0.
\end{align}
 Using (\ref{eq-masstrenpluripolar}), we can consider $w$ as a bounded function with respect to the trace measure of $\langle (\ddc \varphi_j+\omega) \dot{\wedge} T \rangle$. 
Using (\ref{inclu-psikkhacko}), we have 
\begin{align*} 
w \psi_k \ddc \varphi_{jk} \wedge T &= w\bold{1}_{\{\varphi_j >-k\}} \psi_k \langle \ddc \varphi_{j} \dot{\wedge} T \rangle\\
&=w\langle \ddc \varphi_{j} \dot{\wedge} T \rangle+ w(\bold{1}_{\{\varphi_j >-k\}}\psi_k-1 ) \langle \ddc \varphi_{j} \dot{\wedge} T \rangle.
\end{align*} 
The second term in the right-hand side of the last equality converges weakly to $0$ as $k \to \infty$ by the fact that  $\psi_k \to 1$ pointwise outside $A$  as $k \to \infty$ and  Lebesgue's dominated convergence theorem. Hence 
$$ w\langle \ddc \varphi_{j} \dot{\wedge} T \rangle =  \lim_{k \to \infty} w \psi_k \ddc \varphi_{jk} \wedge T.$$
Applying the last equality to $j=1,2$, and using $v= \varphi_1- \varphi_2$, we obtain
\begin{align} \label{eq-psikddcvarphi} 
w \langle \ddc v \dot{\wedge} T\rangle= \lim_{k \to \infty} w \psi_k \ddc v_k \wedge T=\lim_{k \to \infty} w_k \psi_k \ddc v_k \wedge T.
 \end{align}
 Here in the second equality we used the fact that $w=w_k$ on $\{\varphi_3>- k\} \cap \{\varphi_4>-k\}$ which contains $\{\psi_k \not =0\}$.   We also have an analogous formula by exchanging the roles of $v,w$.  Thus,
 \begin{align}\label{eq-psikddcvarphi2} 
w \langle \ddc v \dot{\wedge} T\rangle - v \langle \ddc w \dot{\wedge} T\rangle=  \lim_{k \to \infty} \psi_k (w_k \ddc v_{k}- v_k \ddc w_k) \wedge T.
 \end{align}
By integration by parts for bounded psh functions, we have 
 \begin{multline}\label{eq-interspikaddcv}
\int_X \psi_k (w_k \ddc v_k- v_k \ddc w_k) \wedge T =-\int_X w_k d \psi_k \wedge \dc v_k \wedge T +\\ \int_X v_k d \psi_k \wedge \dc w_k \wedge T. 
 \end{multline}
 Denote by $I_1, I_2$ the first and second term in the right-hand side of the last equality. We will check that $I_j \to 0$ as $k \to \infty$ for $j=1,2$.  Using the  Cauchy-Schwarz inequality, the boundedness of  $v_k,w_k$ and Lemma \ref{le-dvdcv}, we infer
\begin{align*}
|I_1|&\le  \bigg(\int_X  d\psi_k \wedge  \dc\psi_k \wedge  T\bigg)^{\frac{1}{2}} \times  \bigg(\int_X   |w_k|^2 d v_k \wedge \dc v_k \wedge  T\bigg)^{\frac{1}{2}} \\
& \lesssim \bigg(\int_X d\psi_k \wedge  \dc\psi_k \wedge  T\bigg)^{\frac{1}{2}}.
\end{align*}
Recall that  $\{\lim_{k \to\infty} \psi_k <1\}$ is equal to the complete pluripolar set $\{\psi= -\infty\}$. Using this, Theorem \ref{the-th-increasingsequenceMa} and the fact that $T$ has no mass on $\{\psi= -\infty\}$,  we get
$$\lim_{k \to \infty} d\psi_k \wedge  \dc\psi_k \wedge T=\lim_{k \to \infty} (\ddc \psi^2_k - \psi_k \ddc \psi_k) \wedge  T=0$$
Thus we obtain
\begin{align} \label{eq-loimdcspik}
\lim_{k \to \infty } I_1 =0.
\end{align} 
By similarity, we also get $I_2 \to 0$ as $k \to \infty$. Combining this with (\ref{eq-interspikaddcv}) and (\ref{eq-psikddcvarphi})  gives the first desired equality of (\ref{eq-intebyparts}).  We prove the second one similarly as follows. Put $u:= v+w$, and 
$$u_k:= \max\{\varphi_1+ \varphi_3, -k\} - \max\{\varphi_2+ \varphi_4, -k\}.$$
By (\ref{inclu-psikkhacko}) observe that 
$$\bold{1}_{\{\psi_k >0\}} \max\{\varphi_1+ \varphi_3, -2k\} = \bold{1}_{\{\psi_k >0\}}(\varphi_{1k}+ \varphi_{3k})$$
and a similar equality for $\varphi_2, \varphi_4$ also holds. Thus, by plurifine locality, we get 
\begin{align*}
2 \langle d v \wedge \dc w \dot{\wedge} T\rangle &= \langle d u \wedge \dc u \dot{\wedge} T\rangle- \langle d v \wedge \dc v \dot{\wedge} T\rangle-\langle d w \wedge \dc w \dot{\wedge} T\rangle \\
&= \lim_{k \to \infty} \psi_k\big(\langle d u_{2k} \wedge \dc u_{2k} \dot{\wedge} T\rangle- \langle d v_{k} \wedge \dc v_{k} \dot{\wedge} T\rangle -\\
& \quad  \langle d w_{k} \wedge \dc w_{k} \dot{\wedge} T\rangle \big)\\
&= \lim_{k \to \infty} \psi_k\big(\langle d (v_{k}+w_k) \wedge \dc (v_{k}+w_k) \dot{\wedge} T\rangle- \langle d v_{k} \wedge \dc v_{k} \dot{\wedge} T\rangle -\\
& \quad  \langle d w_{k} \wedge \dc w_{k} \dot{\wedge} T\rangle \big).
\end{align*}
Consequently
\begin{align*}
\langle d v \wedge \dc w \dot{\wedge} T\rangle = \lim_{k \to \infty} \psi_k \langle d v_{k} \wedge \dc w_{k} \dot{\wedge} T\rangle.
\end{align*}
It follows that 
 \begin{align*} 
\int_X v \langle \ddc w \dot{\wedge} T\rangle +\langle d v \wedge \dc w \dot{\wedge} T\rangle &=  \lim_{k \to \infty} \int_X  \psi_k (v_k \ddc w_{k}+ d v_k \wedge \dc w_k) \wedge T\\
&= -\lim_{k \to \infty}\int_X v_k d\psi_k\wedge  \dc w_{k} \wedge T
 \end{align*}
 which is equal to $0$ by analogous arguments as in the proof of (\ref{eq-loimdcspik}). This finishes the proof.
\end{proof}

\begin{corollary} \label{cor-integraintro} Let $v,w$ be  bounded dsh functions on $X$.  Then, for every closed smooth form $\Phi$ of right bi-degree,  we have  
\begin{multline}\label{eq-intebypartscor}
\int_X w \langle \ddc v \wedge \bigwedge_{j=1}^m T_j \rangle  \wedge \Phi=\int_X v \langle \ddc w \wedge \bigwedge_{j=1}^m T_j \rangle  \wedge \Phi=\\- \int_X \langle dv \wedge \dc w \wedge \bigwedge_{j=1}^m T_j \rangle  \wedge \Phi.
\end{multline}
\end{corollary}

\proof  By writing $\Phi$ as the difference of two closed positive forms, we can assume that $\Phi$ is positive.  The desired formula is a direct consequence of Theorem \ref{th-integra} applied to $T:= \langle T_1 \wedge \cdots \wedge T_m \rangle \wedge \Phi$.
\endproof

We recall that the first inequality of (\ref{eq-intebypartscor}) was proved in \cite[Theorem 1.2]{Lu-comparison-capacity} and \cite{Xia} when $m=n$ and the cohomology classes of $T_j$'s are big. One should notice a crucial point that the integration by parts formulae obtained in \cite{Lu-comparison-capacity,Xia} contain no term involving $dv \wedge \dc w$. Such a term is essential in applications, especially, in the pluricomplex energy theory. The following result is more general than Theorem \ref{th-integra}. We will need it later.

\begin{theorem} \label{th-integra2chi} Let  $T$ a closed positive current of bi-degree $(n-1,n-1)$ on $X$. Let $v,w$ be bounded $T$-admissible dsh functions on $X$. Let $\chi: \R \to \R$ be a $\cali{C}^3$ function. Then we have  
\begin{align}\label{eq-intebypartschi}
\int_X \chi(w) \langle  \ddc v \dot{\wedge} T\rangle=\int_X v \chi''(w) \langle dw \wedge \dc w \dot{\wedge} T\rangle+\int_X v \chi'(w) \langle \ddc w \dot{\wedge} T\rangle.
\end{align}
\end{theorem}

A quick heuristic reason explaining why (\ref{eq-intebypartschi}) should hold is because $\ddc \chi(w)= \chi''(w) d w\wedge \dc w+ \chi'(w) \ddc w$ if $w$ is a bounded quasi-psh function. 

\proof 
We first note that  \cite[Lemma 5.7]{Viet-generalized-nonpluri} still holds for dsh functions which are the differences of two bounded quasi-psh functions. Now, to obtain the desired equality,  we just follow the proof of Theorem \ref{th-integra} verbatim with $\chi(w)$ in place of $w$. The only thing we need to  clarify is the computation concerning $\ddc \chi(w_k) \wedge T$. To this end, it suffices to use  \cite[Lemma 5.7]{Viet-generalized-nonpluri} because $w_k$ is the difference of two bounded quasi-psh functions. This finishes the proof.   
\endproof

\section{Currents with finite relative energy} \label{sec-weighted}

Let $X$ be a compact K\"ahler manifold.   Let $\alpha_1, \ldots,\alpha_m$ be pseudoeffective $(1,1)$-classes of $X$ and  $T$ a closed positive current on $X$.
Let $P_j$ be a closed positive $(1,1)$-current in the class $\alpha_j$ for $1 \le j \le m$. Put $\mathbf{P}:= (P_1, \ldots, P_m)$. We define $\mathcal{E}_{\mathbf{P}}(T)$ to be the set of $m$-tuple $(T_1, \ldots, T_m)$ of closed positive $(1,1)$-currents such that $T_j\in \alpha_j$ and $T_j$ is more singular than $P_j$ and 
$$\{\langle \bigwedge_{j=1}^m T_j \dot{\wedge} T\rangle\}=\{\langle \bigwedge_{j=1}^m P_j \dot{\wedge} T\rangle\}.$$
Notice that for every  current $P'_j$  in $\alpha_j$ such that $P'_j$ has the same singularities as $P_j$ for $1 \le j \le m$, by the monotonicity of relative non-pluripolar products (see (\ref{ine-mono-relativenopluri})), we have 
$$ \mathcal{E}_{\mathbf{P}}(T)= \mathcal{E}_{\mathbf{P}'}(T).$$ 
Hence, when $P_j$ has minimal singularities in $\alpha_j$ for $1 \le j \le m$, we recover  the class $\mathcal{E}(\alpha_1, \ldots,\alpha_m,T)$  of currents of full mass intersection introduced in \cite{Viet-generalized-nonpluri,BEGZ} because we have 
$$\mathcal{E}_{\mathbf{P}}(T)=\mathcal{E}(\alpha_1, \ldots,\alpha_m,T)$$
in this case.

Let $\chi \in \mathcal{W}^-$. Write $P_j= \ddc \varphi_j+ \theta_j$, where $\theta_j$ is a smooth form and $\varphi_j$ is a negative $\theta_j$-psh function. Let $(T_1, \ldots, T_m) \in \mathcal{E}_{\mathbf{P}}(T).$ Let $u_j$ be a negative $\theta_j$-psh function so that  $T_j= \ddc u_j+ \theta_j$ and $u_j \le \varphi_j$ for $1 \le j\le m$. 
For a negative Borel function $\xi$, we put
\begin{align}\label{eq0-dnenergy}
E_{\xi, \mathbf{P}}(T_1, \ldots,T_m;T):= \sum_{J}\int_X -\xi   \big \langle \bigwedge_{j \in J} T_j  \wedge \bigwedge_{j \not \in J} P_j \dot{\wedge} T \big\rangle,
\end{align}
where the sum is taken over every subset $J$ of $\{1, \ldots, m\}$. \emph{The $(T,\mathbf{P})$-relative joint $\chi$-energy} of $T_1, \ldots, T_m$ is, by definition,  $E_{\xi, \mathbf{P}}(T_1, \ldots,T_m;T)$, where 
$$\xi:= \chi\big((u_1- \varphi_1)+ \cdots+ (u_m- \varphi_m)\big).$$     
The last energy depends on the choice of $u_j,\varphi_j$ but its finiteness does not. That notion generalizes those in \cite{BEGZ,GZ-weighted,Viet-generalized-nonpluri}, see also \cite{Cegrell} for the local setting. 

For every closed positive $(1,1)$-current $P$, let $I_P$ be the set of $x\in X$ so that the potentials of $P$ are equal to $-\infty$ at $x$. Note that $I_P$ is a complete pluripolar set.  By Proposition \ref{pro-sublinearnonpluripolar} $(iv)$, the right-hand side of (\ref{eq0-dnenergy}) remains unchanged if we replace $T$ by $\bold{1}_{X \backslash \bigcup_{j=1}^m I_{P_j}}T$. Hence, in practice, we can assume $T$ has no mass on $\bigcup_{j=1}^m I_{P_j}$.   We denote by  $\mathcal{E}_{\chi, \mathbf{P}}(T)$ the subset of  $\mathcal{E}_{\mathbf{P}}(T)$ containing every $(T_1,\ldots, T_m)$ such that their $(T,\mathbf{P})$-relative joint  $\chi$-energy is finite.

 Here is a monotonicity for the class $\mathcal{E}_{\chi, \mathbf{P}}(T)$ when $P_j=P$ for every $1 \le j \le m$. This generalizes \cite[Theorem 5.8]{Viet-generalized-nonpluri}.

\begin{theorem} \label{th-jointenergy} Let $P=\ddc \varphi+ \theta$ be a closed positive $(1,1)$-current and $\mathbf{P}:=(P, \ldots, P)$ ($m$ times $P$).  Let $\chi \in \mathcal{W}^-$ with $|\chi(0)| \le 1$. Let $(T_{1}, \ldots, T_{m}) \in \mathcal{E}_{\chi, \mathbf{P}}(T)$ and  $(T'_{1}, \ldots, T'_{m}) \in \mathcal{E}_{\mathbf{P}}(T)$  such that $T_j= \ddc u_j + \theta$, $T'_j= \ddc u'_j+ \theta$ such that $u_j, u'_j$ are $\theta$-psh and $u_j \le u'_j \le \varphi$. Put  
$$\xi:= \chi\big((u_1- \varphi)+ \cdots +(u_m- \varphi)\big).$$
Then we have 
$$E_\xi(T'_1, \ldots,T'_m;T) \le  c_1 E_\xi(T_1, \ldots,T_m;T)+ c_2,$$
for some constants $c_1, c_2>0$ independent of $\chi$. In particular, $(T'_{1}, \ldots, T'_{m}) \in \mathcal{E}_{\chi,\mathbf{P}}(T)$.
\end{theorem}

\proof    As mentioned above, we can assume that $T$ has no mass on $I_{P}= \{\varphi= -\infty\}$.    Note here that $\{\varphi = -\infty\} \subset \{u_j= -\infty\}$.  Put
$$u_{jk}:= \max\{u_j, \varphi-k\}- \varphi$$
which is a  bounded dsh function and $$T_{jk}:= \ddc u_{jk}+ P.$$ Observe that $u_{jk}$'s are admissible with respect to $T$.  
Define $u'_{jk}, T'_{jk}$ similarly.  Put
 $$v:= \sum_{j=1}^m (u_j- \varphi_j), \quad v_k:= \max\{v, - k \}, \quad \xi_k= \chi(v_k).$$
 Note that $\xi= \chi(v)$.   With these notations and  a suitable integration by parts ready in our hands (Theorem \ref{th-integra2chi})  replacing \cite[Lemma 5.7]{Viet-generalized-nonpluri}), the proof goes exactly as in the proof of \cite[Theorem 5.8]{Viet-generalized-nonpluri}. The only minor modifications are:  the K\"ahler form $\omega$ is substituted by $P$ and  the wedge products appearing in the proof of  \cite[Theorem 5.8]{Viet-generalized-nonpluri} need to be replaced by  $T$-relative non-pluripolar products. This finishes the proof.
\endproof

The following is a direct consequence of Theorem \ref{th-jointenergy}. 

\begin{corollary} \label{cor-bangnnyuaEPPphay} Let $P, \mathbf{P}$ be as in Theorem \ref{th-jointenergy}.  Let $P'$ be a current in $\{P\}$ which is of the same singularity type as $P$. Then, for $\mathbf{P}'=(P', \ldots, P')$ ($m$ times $P'$), we have  
$$\mathcal{E}_{\chi, \mathbf{P}'}(T)=\mathcal{E}_{\chi, \mathbf{P}}(T).$$
\end{corollary}

For every closed positive $(1,1)$-current $P$, we define  the class $\mathcal{E}_{m,P}(T)$ (resp. $\mathcal{E}_{\chi,m,P}(T)$) to be the set of currents  $T_1 \in \{P\}$ such that $(T_1,\ldots,T_1)$ belongs to $\mathcal{E}_{\mathbf{P}}(T)$ (resp. $\mathcal{E}_{\chi, \mathbf{P}}(T)$), where $\mathbf{P}=(P,\ldots, P)$ ($m$ times $P$). The last space was introduced in \cite{Lu-Darvas-DiNezza-mono} when $T$ is the constant function equal to $1$. As in the case of the usual class of currents of full mass intersection  (\cite[Proposition 2.2]{GZ-weighted}), notice that 
$$\mathcal{E}_{m,P}(T)= \bigcup_{\chi \in \mathcal{W}^-}\mathcal{E}_{\chi,m,P}(T).$$
Let $\alpha$ be a pseudoeffective $(1,1)$-class.   
 By Corollary \ref{cor-bangnnyuaEPPphay}, we see that the notion of the weighted class $\mathcal{E}_{\chi,m,P}(T)$ makes sense if we replace $P$ by its equivalent class (in terms of singularity type) of $(1,1)$-currents.  Hence, we can  define  $\mathcal{E}_{m}(\alpha,T)$ (resp. $\mathcal{E}_{\chi, m}(\alpha,T)$)  to be the set $\mathcal{E}_{m,P}(T)$ (resp. $\mathcal{E}_{\chi, m,P}(T)$), where $P$ is a current with minimal singularities in $\alpha$.  

\begin{theorem}\label{th-plurifine}
Let $U$ be an open subset in $\C^n$.   Let $T$ be a closed positive current on $U$ and $u_{j}, u'_j$ bounded psh functions on $U$ for $1 \le j \le m$, where $m \in \N$. Let $v_j,v'_j$ be  psh functions on $U$ for $1 \le j \le q$. Assume that $u_{j}= u'_{j}$ on $W:=\bigcap_{j=1}^q\{v_j> v'_j\}$ for $1 \le j \le m$. Then we have 
\begin{align}\label{eq-MAequaonstrongplurifine}
\bold{1}_{W} \ddc u_{1} \wedge \cdots \wedge \ddc u_{m} \wedge T=\bold{1}_{W} \ddc u'_{1} \wedge \cdots \wedge  \ddc u'_{m} \wedge T.
\end{align}  
\end{theorem}

\proof If $v_j,v'_j$ are all bounded, then the desired assertion is Theorem 2.9 in \cite{Viet-generalized-nonpluri}. In general, observe that 
$$\{v_j > v'_j\}= \bigcup_{k=1}^\infty \{v_{jk} > v'_{jk}\},$$
where $v_{jk}:= \max\{v_j, -k\}$ and similarly for $v'_{jk}$.  Let $W_k:= \bigcap_{j=1}^q \{v_{jk} > v'_{jk}\}$. We have $W= \bigcup_{k=1}^\infty W_k$ and  $u_j= u'_j$ on $W_k$. Applying  \cite[Theorem 2.9]{Viet-generalized-nonpluri} to $u_j, u'_j,W_k$ gives 
$$\bold{1}_{W_k} \ddc u_{1} \wedge \cdots \wedge \ddc u_{m} \wedge T=\bold{1}_{W_k} \ddc u'_{1} \wedge \cdots \wedge  \ddc u'_{m} \wedge T$$
 for every $k$. Hence, the desired assertion follows. This finishes the proof. 
\endproof

Now, using Theorems \ref{th-jointenergy} and \ref{th-plurifine} instead of \cite[Theorem 5.8]{Viet-generalized-nonpluri} and \cite[Theorem 2.9]{Viet-generalized-nonpluri} respectively, and following arguments in the proof of \cite[Theorems 5.9 and 5.1]{Viet-generalized-nonpluri}, we immediately obtain the following result. 

\begin{theorem}\label{th-convexrelativeclass}  For $\chi \in \mathcal{W}^-$, the sets  $\mathcal{E}_{\chi,m,P}(T)$ and  $\mathcal{E}_{m,P}(T)$ are convex.
\end{theorem}

Finally, we would like to make the following comment.

\begin{remark} Let $\mathcal{W}^+_M$ be the class of weights introduced in \cite[Page 462]{GZ-weighted}. Using arguments from the proof of \cite[Lemma 3.5]{GZ-weighted} and that of Theorem \ref{th-convexrelativeclass}, we can prove the convexity of $\mathcal{E}_{\chi,m,P}(T)$ for $\chi \in \mathcal{W}^+_M$. 
\end{remark}

\subsection*{Acknowledgements}
The author would like to thank referees for their remarks improving the presentation of the paper.

\bibliography{biblio_family_MA,biblio_Viet_papers}
\bibliographystyle{siam}

\bigskip

\end{document}